\documentclass{doublecol-new}

\usepackage[numbers]{natbib}

\usepackage[dvipsnames,svgnames,x11names]{xcolor}
\usepackage[final,markup=bfit,deletedmarkup=xout]{changes}
\definechangesauthor[color=BrickRed]{GBL}

\usepackage{stfloats}
\usepackage{mathrsfs}
\usepackage{graphicx}
\usepackage{subfigure}
\usepackage{multirow}

\theoremstyle{TH}{

}

\theoremstyle{THrm}{

}

\theoremstyle{THhit}{

}

\makeatletter
\def\theequation{\arabic{equation}}

\makeatother

\begin{document}%

\setcounter{page}{1}

\LRH{G.M. Platt et~al.}

\RRH{Basins of attraction and critical curves for Newton-type methods in a phase equilibrium problem}

\VOL{x}

\ISSUE{x}

\PUBYEAR{xxxx}

\BottomCatch


\PUBYEAR{201X}

\subtitle{}

\title{Basins of attraction and critical curves for Newton-type methods in a phase equilibrium problem}

\authorA{Gustavo Mendes Platt $^{\ast}$}
\affA{Grupo de Engenharia e Otimiza\c{c}\~{a}o de Processos Industriais, Escola de Qu\'imica e Alimentos, Universidade Federal do Rio Grande \\ Santo Ant\^{o}nio da Patrulha, RS, Brazil \\
E-mail: gmplatt@furg.br \\ $^{\ast}$ Corresponding author}
\authorB{Fran S\'{e}rgio Lobato}
\affB{Laborat\'{o}rio de Modelagem, Simula\c{c}\~{a}o, Controle e Otimiza\c{c}\~{a}o de Processos Qu\'{i}micos, Faculdade de Engenharia Qu\'{i}mica, Universidade Federal de Uberl\^{a}ndia \\Uberl\^{a}ndia, MG, Brazil \\ E-mail: fslobato@ufu.br}

\authorC{Gustavo Barbosa Libotte and Francisco Duarte Moura Neto}

\affC{Instituto Polit\'{e}cnico, Universidade do Estado do Rio de Janeiro \\ Nova Friburgo, RJ, Brazil \\ E-mail: gustavolibotte@iprj.uerj.br \\ E-mail: fdmouraneto@gmail.com}

\begin{abstract}
Many engineering problems are described by systems of nonlinear equations, which may exhibit multiple solutions, in a challenging situation for root-finding algorithms. The existence of several solutions may give rise to complex basins of attraction for the solutions in the algorithms, with severe influence in their convergence behavior. In this work, we explore the relationship of the basins of attractions with the critical curves (the locus of the singular points of the Jacobian of the system of equations) in a phase equilibrium problem in the plane with two solutions, namely the calculation of a double azeotrope in a binary mixture. The results indicate that the conjoint use of the basins of attraction and critical curves can be a useful tool to select the most suitable algorithm for a specific problem.
\end{abstract}

\KEYWORD{Newton's methods ; Basins of attraction ; Nonlinear systems ; Phase equilibrium.}

\maketitle

\section{Introduction}

Root-finding methods have been extensively used in many engineering fields in the last decades. In many situations, the nonlinear systems to be solved show more than one solution; see\deleted{, for instance,}  \citet{Industrial2018}, \citet{Canadian2015} and \citet{BMSE2016}. In such cases, the construction of basins of attraction for the solutions can be useful to determine convergence patterns for the different algorithms \citep{Zotos2017}.

As pointed by \citet{Bischi2011}, the use of basins of attractions \replaced{is a}{are} common tool\deleted{s} in the study of dynamical systems and in iterative numeric procedures arising from the employment of root-finding methods. \citet{Scott2011} detailed the basins of attraction for several algorithms (including Newton's and Halley's methods) in order to evaluate the impact of this diagram in the convergence properties of the investigated methods, using a polynomial equation with one real root and two complex conjugate roots. \citet{Schneebeli2011} proposed an adaptive Newton method based on a dynamic approach to solve nonlinear systems of equations. These authors also analyzed basins of attractions for the original and modified Newton's methods, illustrating the relationship of the critical curve (where the Jacobian matrix of the system is singular) and some confinement of the basins of attraction. \deleted{Recently,} \replaced{I}{i}n a subsequent work,  \citet{Amrein2014} analized the percentage of convergent iterations and the average convergence order for this modified Newton's method in a nonlinear algebraic system with six roots (in the plane); basins of attraction and the critical set were displayed in the same figures. \added{\citet{Chun2019}} \added{analysed the behavior of several methods of various orders in the calculation of roots of nonlinear equations. They presented the number of function evaluations and also discussed basins of attraction for each algorithm.} \added{\citet{Pet2020}} \added{studied the computational performance of simultaneous methods in root-finding problems, using a computational analysis of the convergence order. They advocate the superiority of these classes of algorithms when compared to ``deflation'' techniques, in which each root is found separately.}

In \replaced{c}{C}hemical \replaced{e}{E}ngineering context, basins of attraction were previously studied for some interesting problems, including calculations in complex domains \citep{Lucia1993}, multiple dew point solutions \citep{Lucia1990}, double retrograde vaporization \citep{Platt2012} and critical point calculations \citep{Parajara2017}. \citet{Green1993} presented the fractal behavior in the stationary point calculation problem. Furthermore,  \citet{Guedes2011} presented an initial analysis of the use of basins of attraction for Newton's methods in the calculation of a double azeotrope in the binary system benzene + hexafluorobenzene, but only as a motivation for the use of metaheuristics for solving nonlinear algebraic systems. \added{In fact, metaheuristics can be employed in the solution of nonlinear algebraic problems, converted to a scalar fitness function, but with a higher computational cost when compared to deterministic approaches} \added{\citep{Platt2018, Cheng2018, Ma2019, Chen2019}.}

In this work, we investigate the basins of attraction for several Newton-type algorithms in the calculation of a double azeotrope in the system 1,1,1,-2,3,4,4,5,5,5-decafluoropentane (HFC-4310 mee) + oxolane (THF), also analyzing the effect of the critical curve -- the locus where the Jacobian of the system of equations is non-invertible -- of the nonlinear system in the basins of attraction. The calculation of azeotropes in a binary mixture (a mixture with two substances) can be represented by a 2 $\times$ 2 algebraic system with the molar fraction of one component and the system temperature as unknowns (under specified pressure). Thus, considering this situation, the convergence analysis for
deterministic algorithms (as Newton's methods) can be conducted with the aid of basins of attraction in a two-dimensional domain.

\added{The proposed methodology is useful considering the development of robust computational frameworks devoted to the solution of phase equilibrium problems, such as azeotropy. The accurate and robust calculation of azeotropic coordinates is applied, for instance, in refrigeration systems} \added{\citep{Zhao2019}.}

\section{Methodology and Modelling of the Problem}

In this section, the Newton's method and some variants are presented in a concise form. We also provide the formulation of the azeotrope calculation problem. The last subsection of this section is dedicated to the description of the computational convergence orders for the studied methods.

\subsection{Newton's Method (CN)}

Newton's method (or classic Newton's method, CN) is probably the best known and most used to find the
roots in a nonlinear equations system $f(X)=0$. Mathematically, the iterative formula that represents it can be developed by using a geometric approach, or by truncating Taylor's formula, resulting in the following iteration
\begin{equation}
X^{k+1}=X^k+\mathit{\delta X}^k
\label{eq:step_newton}
\end{equation}
\begin{equation}
J(X^k)\mathit{\delta X}^k=-f(X^k)
\label{eq:sist_newton}
\end{equation}
where $k$ is the iteration counter, $J$ is the Jacobian matrix and $\delta
X$ is the step size. 

Given an initial estimate $X^0$, and assuming \added{that the} estimate $X^k$ has been computed, in each iteration, it is necessary to solve a linear system to find $\delta X^k$ and $X^{k+1}$. If a predefined stop criterion is reached, this procedure is finalized. Otherwise, the computed value is considered as
a new estimate to update the value of the solution. To solve this linear system, direct or iterative
approaches can be used, such as Jacobi, Gauss-Seidel, Successive Overrelaxation, Accelerated Overrelaxation and Krylov subspace methods \citep{Ramos2015, Amiri2018}\deleted{)}.

The Newton's method converges quadratically when a good initial estimate is used. This represents its
main advantage. On the other hand, Newton's method may fail if the determinant of the Jacobian matrix is either zero or very small. In addition, it may fail also to converge to the desired root in case the initial estimate is far from this root. In some cases the iterations oscillate and the convergence may be very slow near roots of multiplicity greater than one \citep{Liu2015, Ramos2015, Amiri2018}\deleted{)}.

In order to overcome these disadvantages, various modifications in the Newton's method can be found in
the literature. \replaced{The following}{Next} subsections present three variants: Global Newton Method with Residual based Convergence Criterion, Global Newton Method with Error Oriented Convergence Criterion, and Jacobian-free Newton--Krylov Method.

\subsection{Global Newton Method with Residual based Convergence Criterion (NLEQ-RES)}

During the execution of the classical Newton's Method, the initial estimate is an important factor
that influence the convergence process. To minimize this influence, Global Newton methods pursue two approaches, damping or adaptive trust region strategies \citep{Deuflhard2004}. Here, the Global Newton method with residual based convergence criterion (NLEQ-RES) will be briefly presented.

In this approach, a modified problem is solved:
\begin{equation}
X^{k+1}=X^k+\lambda ^k\mathit{\delta X}^k
\label{eq:damped}
\end{equation}
\begin{equation}
J(X^k)\mathit{\delta X}^k=-f(X^k)
\label{eq:damped_2}
\end{equation}
where $\lambda^k$ is the damping factor ($0 < \lambda^k \leq 1$).

The damping factor is updated considering the following procedure \citep{Deuflhard2004}. Initially, an initial estimate $X^0$ and an initial damping factor (with $\lambda^k \leq 1 $) are defined by the user. If the convergence test is not satisfied, the $\delta X^k$ value is updated by using the Equation (\ref{eq:damped_2}). Then, the damping factor is computed considering the adaptive trust region strategy \citep{Deuflhard2004}:
\begin{eqnarray}
\lambda^k=\min(1,\mu^k) \\
\mu ^k=\frac{||f(X^{k-1})||}{||f(X^k)||}\mu
^{k-1}
\end{eqnarray}

Next, to update this value, a regularity test is performed (if $\lambda^k < \lambda^{min}$, the iterative process is finalized). Otherwise, compute a new value for $X^{k+1}$
by using Equation (\ref{eq:damped}).

During the iterative process, the monitoring quantities are computed as (regularity test):
\begin{eqnarray}
\theta^k=\frac{||f(X^{k+1})||}{||f(X^k)||} \\
\mu^k=\frac{0.5||f(X^k)||(\lambda^k)^2}{||f(X^{k+1})-(1-\lambda ^k)f(X^k)||}
\end{eqnarray}

If $\theta^k \geq 1$ (or, if restricted: $\theta^k > 1 - \frac{\lambda^k}{4}$) then replace
$\lambda^k$ by $\lambda^k_p = \min(\mu^k, 0.5 \lambda^k)$ and go to Regularity test. Else $\lambda^k_p = \min(1,\mu^k)$. If $\lambda^k_p = \lambda_k = 1$ and $\theta^k < \theta^{max}$, then the iterative process is interrupted. Else if $\lambda^k_p \geq 4$, replace $\lambda^k$ by $\lambda^k_p$ and
compute a new value for $X^{k+1}$ by using Equation (\ref{eq:damped}). Else accept $X^{k+1}$ as a new iterate and update $k = k + 1$.

\subsection{Global Newton Method with Error Oriented Convergence Criterion (NLEQ-ERR)}

The Global Newton method with error-oriented convergence criterion (NLEQ-ERR) is also an approach that
consists in updating the damping factor. Its procedure is very similar to the previous one. Consider the following modified nonlinear system as function of the damping factor $\lambda^k$:
\begin{equation}
X^{k+1}=X^k+\lambda ^k \delta X^k
\label{eq:newton_NLEQ_ERR}
\end{equation}
\begin{equation}
J(X^k)\mathit{\delta X}^k=-f(X^k)
\label{eq:damped_NLEQ_ERR}
\end{equation}

Initially, as in NLEQ-RES method, both initial estimate $X^0$ and damping factor $\lambda^0$ are defined by the user. If the convergence test is not satisfied, the $\delta X^k$ value is updated by using Equation (\ref{eq:damped_NLEQ_ERR}). The damping factor is computed considering the adaptive trust region strategy \citep{Deuflhard2004} (for $k > 0$):
\begin{eqnarray}
\lambda ^k=\min(1,\mu ^k) \\
\mu ^k=\frac{||\delta X^{k-1}||||\delta X^k||}{||\delta X^k- \delta X^k||||\delta X^k||}\lambda^{k-1}
\end{eqnarray}
with
\begin{equation}
J(X^k)\delta X^{k+1}=-f(X^{k+1})
\end{equation}

Next, to update this value, a regularity test is performed (if $\lambda^k < \lambda_{min}$ , the iterative process is finalized). Otherwise, compute a new value for $X^{k+1}$ by Equation (\ref{eq:newton_NLEQ_ERR}).

During the iterative process, the monitoring quantities are computed as (regularity test):
\begin{eqnarray}
\theta ^k=\frac{||\delta X^{k+1}||}{||\delta X^k||} \\
\mu ^k=\frac{0.5||\delta X^k||(\lambda ^k)^2}{||\delta X^{k+1}-(1-\lambda^k)\delta X^k||}
\end{eqnarray}

If $\theta ^k \geq 1$ (or, if restricted: $\theta^k > 1 - \frac{\lambda^k}{4}$) then replace
$\lambda^k$ by $\lambda^k_p = \min(\mu^k, 0.5\lambda^k)$ and go to regularity test. Else $\lambda^k_p = \min(1,\mu^k)$. If $\lambda^k_p = \lambda^k = 1$ and $|| \delta X^{k+1}|| <\varepsilon $ ($\varepsilon$ is a small number and, in our implementation, equals to $1 \times 10^{-12}$), then the iterative process is interrupted and the solution is found. If $\lambda^k_p \geq 4$, replace $\lambda^k$ by $\lambda^k_p$ and compute a new value for $X^{k+1}$ by using the Equation (\ref{eq:newton_NLEQ_ERR}). Else accept $X^{k+1}$ as new iterate and update $k = k + 1$.

\subsection{Jacobian-free Newton--Krylov Method (JFNK)}

The Jacobian-free Newton-Krylov Method (JFNK) is a variant of classical Newton Method. In general,
this approach consists in solving a series of linearized Newton correction equation (Equation (\ref{eq:sist_newton})) \citep{Knoll2004}. In the JFNK strategy, the linear system, Equation (\ref{eq:sist_newton}), could be effectively solved with a Krylov's method (KM). In this approach, only a matrix-vector product is required and can be approximated as \citep{Knoll2004}:
\begin{equation}
J^kv=\frac{f(X^k+\varepsilon v)-f(X^k)}{\varepsilon}
\end{equation}
where $v$ is the Krylov vector and $\varepsilon $ is a small number (again equal to $1 \times 10^{-12}$ in our implementation). The KM method approximates the linear solution in a subspace of the form
$(r_0; Jr_0; J^2r_0, \ldots)$  where $r_0$ is the initial linear residual with an initial guess $\delta X^k$ usually taken to be zero. Because Krylov's method requires only the product of the system matrix (the Jacobian) and a vector, the Jacobian matrix does not need to be formed and stored explicitly \citep{Lemieux2011,Zou2016}.

After this correction, $\delta X^k$ is obtained by solving the linear system (Equation (\ref{eq:sist_newton})), and a new position ($k+1$) is updated according to Equation (\ref{eq:step_newton}).

\subsection{Method of \citet{Babajee2012} (BA) }

Next, some higher-order Newton-type methods will be presented, in order to evaluate the relation of
the order of convergence of such methods with the convergence rate and the computational cost related to the construction of the basins of attraction. Among these methods, the first one to be approached in this work will be the method proposed by \citet{Babajee2012}. In this two-step iterative scheme, the main favorable feature of sequence convergence is that it is not required to evaluate the second or higher order derivatives, such that only by evaluating two first-order derivatives at each step, fourth-order convergence can be achieved. Therefore, the modification proposed by \citet{Babajee2012} to Newton's method is calculated through the relation
\begin{equation}
X^{k+1}=X^k-W\left(X^k\right)A_1\left(X^k\right)^{-1}f\left(X^k\right)
\end{equation}
in such a way that the corrector step is calculated through
\begin{eqnarray}
A_1\left(X^k\right)=\frac 1
2\left(J\left(X^k\right)+J\left(Y\left(X^k\right)\right)\right) \\
W\left(X^k\right)=I-\frac 1 4\left(\tau \left(X^k\right)-I\right)+\frac 3 4\left(\tau
\left(X^k\right)-I\right)^2
\end{eqnarray}

In turn, $I$ represents the $n\times n$ identity matrix and the terms $\tau \left(X^k\right)$ and $Y\left(X^k\right)$ are given by
\begin{eqnarray}
\tau
\left(X^k\right)=J\left(X^k\right)^{-1}J\left(Y\left(X^k\right)\right) \\
Y\left(X^k\right)=X^k-\frac 2 3J\left(X^k\right)^{-1}f\left(X^k\right)
\end{eqnarray}

\subsection{Method of \citet{Grau2011} (GS)}

Through expansions in formal developments in power series of the function, \citet{Grau2011} proposed a fifth-order method, performing the iterative process through the correction factor calculated by
\begin{equation}
X^{k+1}=Z\left(X^k\right)-J\left(Y\left(X^k\right)\right)^{-1}f\left(Z\left(X^k\right)\right)
\end{equation}
with $Z\left(X^k\right)$ and $Y\left(X^k\right)$ given by

\begin{eqnarray}
Z\left(X^k\right) = \nonumber \\ 
X^k-\frac {1}{2}\left(J\left(X^k\right)^{-1}+J\left(Y\left(X^k\right)\right)^{-1}\right)f\left(X^k\right) \\
Y\left(X^k\right)=X^k-J\left(X^k\right)^{-1}f\left(X^k\right)
\end{eqnarray}

\subsection{Method of \citet{Cordero2012} (CO)}

\citet{Cordero2012} proposed a sixth-order method where the correction factor in the iterative process is computed by following relation
\begin{eqnarray}
X^{k+1} =\nonumber \\
Z\left(X^k\right)-(J\left(X^k\right)-2J(Y(X^k)))^{-1}f\left(Z(X^k)\right)))
\end{eqnarray}
where $Z\left(X^k\right)$ and $Y\left(X^k\right)$ are given by
\begin{eqnarray}
Z\left(X^k\right)=X^k-\left(J\left(X^k\right)^{-1}+2J\left(Y\left(X^k\right)\right)^{-1}\right)\times \nonumber \\
(3F\left(X^k\right)-4F\left(Y\left(X^k\right)\right)) \\
Y\left(X^k\right)=X^k-\frac 1 2J\left(X^k\right)^{-1}f\left(X^k\right)
\end{eqnarray}

\subsection{Method of \citet{Madhu2017} (O4N, O5N and O6N) }

Recently, \citet{Madhu2017} proposed new higher order iterative methods to solve a  system of nonlinear equations where any information about second derivatives is necessary. The fourth order method (O4N) is formulated as
\begin{equation}
X^{k+1}=\mathit{G4th}\left(X^k\right)=X^k-H\left(X^k\right)(J(Y\left(X^k\right)))
\end{equation}
where $H\left(X^k\right)$ is given by
\begin{equation}
 H\left(X^k\right)=I-\frac 3 4\left(\tau \left(X^k\right)-I\right)+\frac 9 8(\tau
\left(X^k\right)-I)^2
\end{equation}
and
\begin{eqnarray}
\tau \left(X^k\right)=J\left(X^k\right)^{-1}J\left(Y\left(X^k\right)\right) \\
Y\left(X^k\right)=X^k-\frac 2 3J\left(X^k\right)^{-1}f\left(X^k\right)
\end{eqnarray}

The fifth order method (O5N) is formulated as
\begin{equation}
X^{k+1}=\mathit{G5th}\left(X^k\right)=\mathit{G4th}\left(X^k\right)-(J\left(X^k\right))
\end{equation}

Finally, the sixth order method (O6N) is formulated as
\begin{eqnarray}
X^{k+1}=\nonumber \\
\mathit{G6th}\left(X^k\right)=\mathit{G4th}\left(X^k\right)-T\left(X^k\right)(J\left(X^k\right))
\end{eqnarray}
where
\begin{equation}
T\left(X^k\right)=I-\frac 3 2\left(\tau \left(X^k\right)-I\right)+\frac 1 2(\tau
\left(X^k\right)-I)^2
\end{equation}

\subsection{Some geometric notions on convergence}

Given a nonlinear function $f:\Lambda {\subseteq}R^n\rightarrow R^n$, it is possible that the nonlinear equation $f\left(X\right)=0$ has multiple solutions, say, $X_1^{\ast },{\dots},X_m^{\ast }$. Assuming that one uses Newton's method to search for the solutions, one can write, from Equations (\ref{eq:step_newton}) and (\ref{eq:sist_newton}), $X^{k+1}=X^k-\left[J\left(X^k\right)\right]^{-1}f(X^k)$. We want to have an idea of what can happen with this iterative method with respect to converging to the solutions.

Let $N\left(X\right)=X-\left[J\left(X\right)\right]^{-1}f(X)$. The Jacobian $J\left(X\right)$ is noninvertible if and only if $\det\left(J\left(X\right)\right)=0$, which defines the set of critical curves of $f$, $C^0(f)$, a closed set in the domain of $f$. Newton's iteration can be rewritten as $X^{k+1}=N\left(X^k\right)$. The sequence thus defined, $X^0, X^1, X^2, \ldots$ is called the orbit of $X^0$ under the iteration defined by the function $N$. However, the orbit is not always defined for all $k$ since it might happen that $X^k{\in}C^0\left(f\right)$, for some $k$, even $k=0$ (and the iteration cannot proceed to compute the next term $X^{k+1}$); in that case we say that  $X^0{\in}C^k\left(f\right)$, which, still, is a closed set. Denote by $C\left(f\right)={\bigcup}_{k=0}^{{\infty}}C^k\left(f\right)$ the set of $X^0$ whose orbit hits $C^0\left(f\right)$ for a certain value of the iteration counter. Now, the set where Newton's iteration is defined is $\Omega =\Lambda \backslash C\left(f\right)$, thus the domain of $f$ splits in two disjoint sets, $\Lambda =\Omega
{\bigcup}^dC\left(f\right)$, one where the iteration procedure can be continued \textit{ad infinitum}, and the other where it has to be halted.

One says that $X^0{\in}\Omega $ is in the basin of attraction of a solution of the nonlinear equation $f\left(X\right)=0$, by the iterative method defined by $N$ if its orbit converges to a point $\underline X$. Then necessarily $\underline X$ is a fixed point of $N$, $\underline X=N(\underline{X})$, and hence it is a solution of the nonlinear equation, $f\left(\underline X\right)=0$. Denote by $B_{X_i^{\ast }}(N)$ the basin of attraction of $X_i^{\ast }$. It may happen that $B_{X_i^{\ast
}}(N)$ is empty. Moreover, if $X_i^{\ast }{\neq}X_j^{\ast }$ then $B_{X_i^{\ast }}\left(N\right){\cap}B_{X_j^{\ast }}\left(N\right)={\emptyset}$. Some transitions from a basin of attraction of one point to another may be traversing the critical set of $f$.

Another set in this respect is the Fatou set of the iteration defined by $N$, $F(N)$, which is the set of points that behave similarly as the iteration proceed, i.e., their orbits remain close, or go to infinity. More precisely, a point  $X^0{\in}\Omega $ is in the Fatou set of the iteration, $X^0{\in} F(N)$, if given $\varepsilon >0$, there is a  $\delta >0$ such that for any sufficiently close point  $Y^0{\in}\Omega $, their orbits are close
\begin{eqnarray}
\|Y^0-X^0\|<\delta \Rightarrow \|Y^k-X^k\|= \nonumber \\
\|N^{(k)}(Y^0)-N^{(k)}(X^0)\|<\varepsilon
\end{eqnarray}
or they go to infinity, $\|Y^0-X^0\|<\delta \Rightarrow \|Y^k=N^{(k)}(Y^0)\|>\varepsilon ^{-1}$
and $\|X^k=N^{(k)}(X^0)\|>\varepsilon ^{-1}$. By construction, Fatou set is
open. The complementary set, a closed set is called Julia set, $J(N)$. The basins of attraction are subsets of the Fatou set,  ${\cup}_iB_{X_i^{\ast }}(N){\subset} F(N)$.

\subsection{Modelling of the problem}

Considering the description presented in the last subsection, we will construct and analyze the basins of attraction and the critical curves in a double azeotrope calculation problem at low pressures. An azeotrope is a thermodynamic condition where a boiling liquid produces a vapor with the same composition of the liquid (obviously, pure fluids attend this condition, and the azeotropy phenomenon refers to that mixtures that follow this unusual behavior) \citep{Swietoslawski1964}. The term azeotrope is derived from Greek and means ``does not change when boiling''. Some particular binary mixtures can exhibit azeotropes and, in rare occasions, more than one azeotrope. In these cases, we obtain the double azeotropy phenomenon (referring to the existence of two azeotropes in a binary mixture; naturally, multicomponent mixtures can show many binary, ternary, quaternary azeotropes).

We will consider, as an example for the azeotrope calculation, the binary mixture formed by 1,1,1,2,3,4,4,5,5,5-decafluoropentane (HFC-4310 mee) + oxolane (THF). As pointed by \citet{Segura2005}, this mixture exhibits two azeotropes at 35 kPa. The coordinates for the azeotropic temperature (Kelvin) and molar fraction of HFC-4310 mee in the liquid phase ($x_1$) are detached in Table \ref{tab:azeotropes}.

\begin{table}[h]
\centering
\begin{tabular}{l l l}
\hline
 & $T(K)$ & $x_1$ \\
\hline
Azeotrope 1 & 309.45 & 0.0923864 \\
Azeotrope 2 & 309.57 & 0.2552517 \\
\hline
\end{tabular}
\caption{Azeotropic coordinates for the binary system HFC-4310 mee + THF at 35 kPa \citep{Canadian2015}}
\label{tab:azeotropes}
\end{table}

The azeotropic condition is, first of all, a coexistence condition. It means that the equations that assure vapor-liquid coexistence in a binary mixture must be attended:

\begin{equation}
f_i^L=f_i^V, \hspace{0.1in} i=1,2,\ldots,c
\end{equation}
where  $f_i^{\alpha }$ \ is the fugacity of component \textit{i} in phase  $\alpha $. Considering an ideal vapor phase (due to the low pressure in the system) and a nonideal liquid phase, the vapor-liquid coexistence can be represented as

\begin{equation}
x_i\gamma _iP_i^{\mathit{sat}}=Py_i, \hspace{0.1in} i=1,2,\ldots,c
\label{eq:raoult}
\end{equation}

In this last equation,  $P_i^{\mathit{sat}}$ is the saturation pressure of a pure component, \textit{P} is the system pressure (35 kPa in this example) and  $\gamma _i$ \ is the activity coefficient of component \textit{i} (that describes the nonideal behavior of the liquid phase; in a general way,  $\gamma _i$= $\gamma_i\left(x_1,x_2,T\right)$). The compositions of liquid and vapor phases are represented by molar fractions,  $x_i$ and  $y_i$, respectively.

Furthermore, since the vapor and the liquid have the same composition, we obtain:
\begin{equation}
x_i=y_i
\label{eq:azeotrope}
\end{equation}

On the other hand, in a binary mixture the molar fractions must obey  $x_1+x_2=1$ \ and 
$y_1+y_2=1$. Thus, the set of Equations (\ref{eq:raoult}) can be written as:
\begin{equation}
x_1\gamma _1P_1^{\mathit{sat}}=Py_1
\label{eq:phase_eq_1}
\end{equation}
\begin{equation}
\left(1-x_1\right)\gamma _2P_2^{\mathit{sat}}=P\left(1-y_1\right)
\label{eq:phase_eq_2}
\end{equation}

With Equation (\ref{eq:azeotrope}) in Equations (\ref{eq:phase_eq_1}) and (\ref{eq:phase_eq_2}), the azeotropy calculation problem is then represented by:

\begin{eqnarray}
\gamma _1P_1^{\mathit{sat}}=P \\
\gamma _2P_2^{\mathit{sat}}=P
\end{eqnarray}

Thus, we have a nonlinear algebraic problem in the plane, with unknowns  $X =\left(x_1,T\right)$\ (under specified pressure).

The activity coefficient (a quantity related to the modeling of the non-idealities in the liquid phase) will be described by a Redlich-Kister model, as pointed by \citet{Segura2005}. The excess molar Gibbs free energy, $G^E$, is represented by:

\begin{eqnarray}
\frac{G^E}{RT}=x_1x_2\left[C_1+C_2\left(x_2-x_1\right)+C_3\left(x_2-x_1\right)^2+ \right. \nonumber \\\left.C_4\left(x_2-x_1\right)^3\right]
\end{eqnarray}
where \textit{R} is the universal gas constant. Parameters  $C_1$ to  $C_4$ are functions of
temperature and are tabulated \cite{Segura2005}. The activity coefficients can be obtained by differentiation of  $G^E$ with respect to the number of mols of each component, as follows:

\begin{equation}
\ln \left(\gamma
_i\right)=\left[\frac{{\partial}\left(\frac{nG^E}{\mathit{RT}}\right)}{{\partial}n_i}\right]_{T,P,n_{j{\neq}i}}
\end{equation}
where  $n=n_1+n_2$. Evaluation of the saturation pressures for pure fluids are obtained using the Antoine equation, with the coefficients presented by \citet{Segura2005}.

Clearly, we are dealing with a nonlinear system in the plane, with  $x_1$ \ and \textit{T} as coordinates. Thus, we can present some useful questions for several algorithms applied to this phase equilibrium problem: (i) How \added{does} the algorithm behave\deleted{s} for a specific initial estimate?; (ii) Is there a clear relationship between the critical curve and the basin of attraction for a specific algorithm?

\subsection{Computational Analysis of the Order of Convergence of the Algorithms}

Since the set of Newton-type methods used in this work to solve the double azeotrope problem and to generate the basins of attraction have different orders of convergence (for more details, see the bibliography of each technique), there was a natural tendency to the relationship between the respective orders of convergence of the algorithms and the characteristics inherent to each basin of attraction. For these type of methods, two approaches are widely used in the literature: the computational order of convergence \citep{Nedzhibov2008} and the approximated computational order of convergence \citep{Grau2010}.

It is known that the analysis of the order of convergence of Newton-type algorithms is conditioned by the fact that the initial estimate must be sufficiently close to the solution. Therefore, given a sequence of points, produced iteratively by a given method, and denoted by  $\{X_i^k\}$, the \replaced{computational}{computacional} order of convergence  $\rho (i,k+1)$  is calculated by

\begin{equation}
\rho \left(i,k+1\right)=\frac{\ln \left|\left(X_i^{k+1}-\alpha _i\right)/\left(X_i^k-\alpha _i\right)\right|}{\ln
\left|\left(X_i^k-\alpha _i\right)/\left(X_i^{k-1}-\alpha _i\right)\right|}
\end{equation}
for each of the variables  $i=1,{\dots},n$. In this approach, one must consider a known solution, given by  $\alpha $. Therefore, the computational order of convergence, calculated by the ratio that relates three consecutive points of the sequence, to multidimensional problems, is given by the average value of  $\rho (k+1)$ \ for each of the variables ($\rho_{avg}$).

When the exact value of the solution  $\alpha $ is not known, the computational order of convergence can be estimated through the rate that considers only consecutive solutions of the sequence $\{X_i^k\}$. In these cases, the approximated computational order of convergence is obtained, which is given by

\begin{equation}
\widehat {\rho }\left(i,k+1\right)=\frac{\ln
\left|\left(X_i^{k+1}-X_i^k\right)/\left(X_i^k-X_i^{k-1}\right)\right|}{\ln
\left|\left(X_i^k-X_i^{k-1}\right)/\left(X_i^{k-1}-X_i^{k-2}\right)\right|}
\end{equation}

In all cases considered, the iteration counter  $k$ was taken as  $p-1$, for a sequence composed of at most $p$ elements and that has at least four elements. Similarly to the previous case, the approximated computational order of convergence of problems with several variables is given by the average between the values calculated for each of the variables $i$ and will be represented as $\widehat {\rho }_{avg}$.

\section{Results and Discussions}

In this section we present the computational experiments in the calculation of the two azeotropes in the binary system HFC-4130 mee + THF at 35 kPa, as represented in Table 1. All methods are evaluated in a grid of initial estimates in the range  $1\times 10^{-6}{\leq}x_1{\leq}1-1\times 10^{-6}$ \ and  $280{\leq}T{\leq}400$, with 40,000 points (a $200\times 200$ grid). If the algorithm converges to the Azeotrope 1, a red point is inserted to identify its basin of attraction. Conversely, when the method converges to Azeotrope 2, a blue point is marked. The yellow portion of the diagrams represents singularities in the Jacobian matrix (a divergence situation, preventing the application of the recurrence relations). The green points indicate the initial conditions that converge to infinity. Finally, when the algorithm shows an oscillatory behavior (without convergence), an orange point is inserted in the figure.

Furthermore, the critical curves are represented by a black continuous line in each figure. In fact, this problem exhibits more than one critical curve, as detailed by \citet{Canadian2015}. Here, we are interested only in portions of the critical curve contained in the feasible domain (i.e., with molar fractions in the interval between zero and one); \citet{Canadian2015} employed a numerical method of inversion of functions to obtain the two azeotropes in the system and, in this case, even non-physical portions of domain were analyzed.

\added{It is important to highlight that the proposed methodology --- the generation of the critical curves and the the construction of the basins of attraction --- is expensive from a computational point of view. In real engineering problems, the calculation of the critical curves is not usually conducted. On the other hand, as we will demonstrate, some convergence patterns of root-finding algorithms can be understood in the light of these computational tools.}

\begin{figure*}[t]
 \centering
 \framebox[2.9in][c]{
 \includegraphics[width=2.7in]{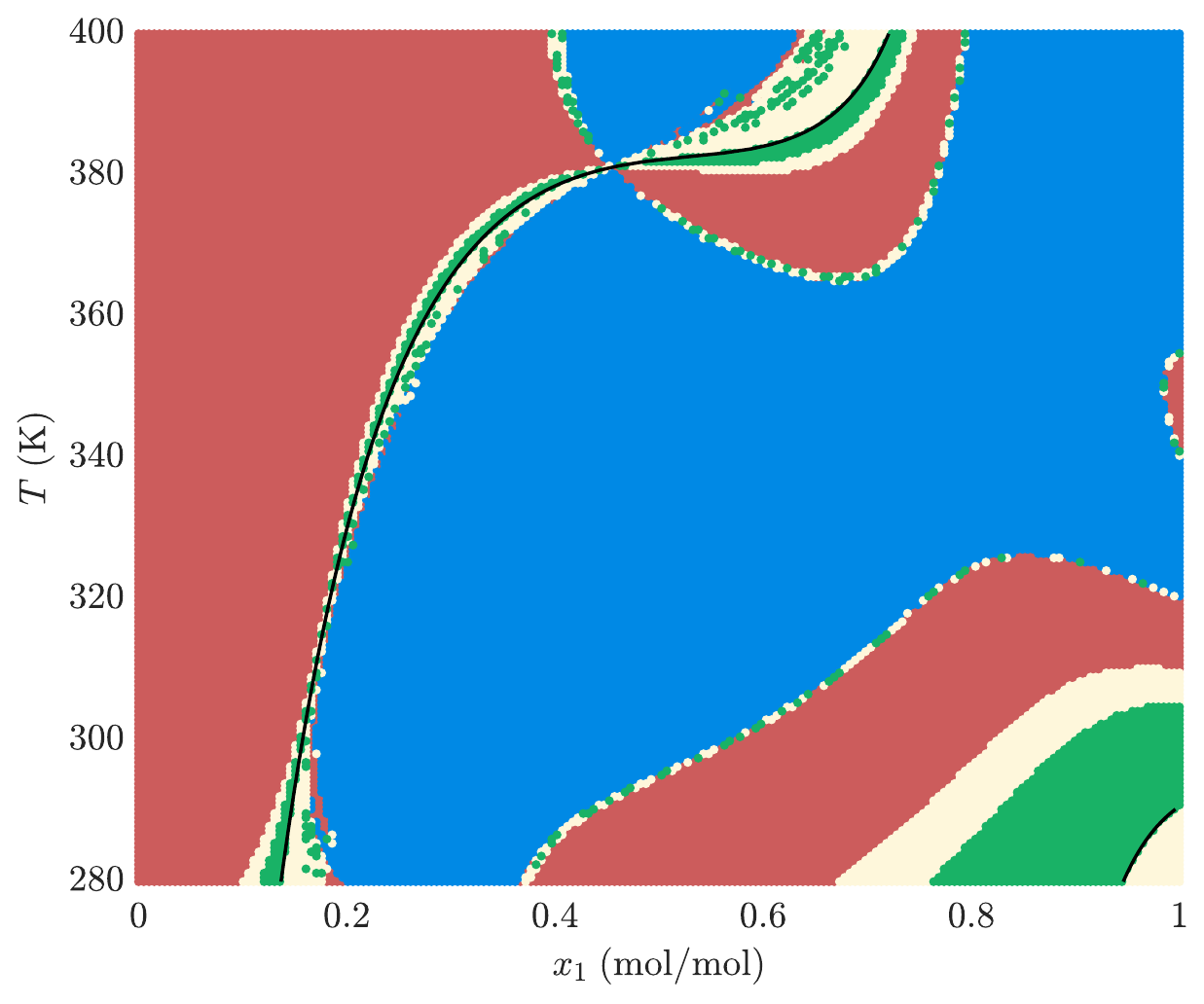}
}
  \caption{Basin of attraction for the classical Newton's method with 200 $\times$ 200 points. Red
points: Azeotrope 1. Blue points: Azeotrope 2. Yellow points: singular Jacobian matrix. Green points: convergence to infinity. Orange points: oscillatory behavior. Black continuous line: critical curve.}
\label{fig:newton}
\end{figure*}

Figure \ref{fig:newton} contains the basin of attraction for the classical Newton method with Jacobian calculated numerically with a five-point formula. We noted that the critical curve essentially ``splits'' the basin of attraction in two big regions: red and blue. Moreover, initial estimates in the neighborhood of the critical curve show undesirable patterns in many occasions (singular Jacobian matrix and convergence to infinity), which is expected since the critical curve represents a fail in the Newton algorithm. A large portion of non-convergence region is found close to pure HFC-4130 mee, i.e., with molar fractions close to one.

\begin{figure*}[t]
 \centering
 \framebox[2.9in][c]{
 \includegraphics[width=2.7in]{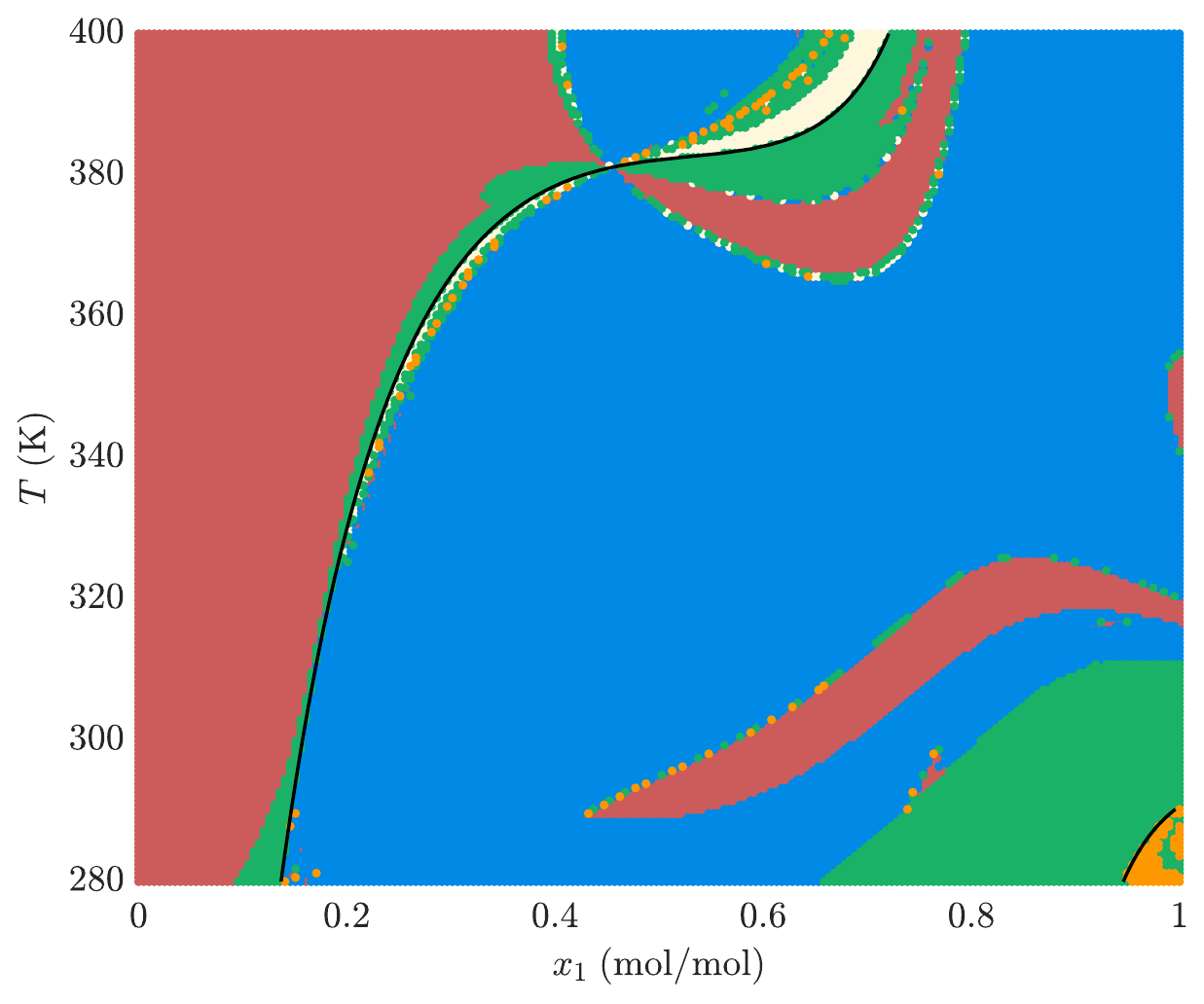}
}
  \caption{Global Newton method with residual based convergence criterion and adaptive trust region
strategy (NLEQ-RES) with 200 $\times$ 200 points. Red points: Azeotrope 1. Blue points: Azeotrope 2. Yellow points:singular Jacobian matrix. Green points: convergence to infinity. Orange points: oscillatory behavior. Black continuous line: critical curve.}
\label{fig:nleq_res}
\end{figure*}

Figure \ref{fig:nleq_res} represents the basin of attraction for algorithm NLEQ-RES. We can note a reduction in the non-convergence region, mainly in the vicinities of the critical curve. On the other hand, the green portion of the figure (initial estimates that promote convergence to infinity) are increased in a comparison with Figure \ref{fig:newton}.

\begin{figure*}[t]
 \centering
 \framebox[2.9in][c]{
 \includegraphics[width=2.7in]{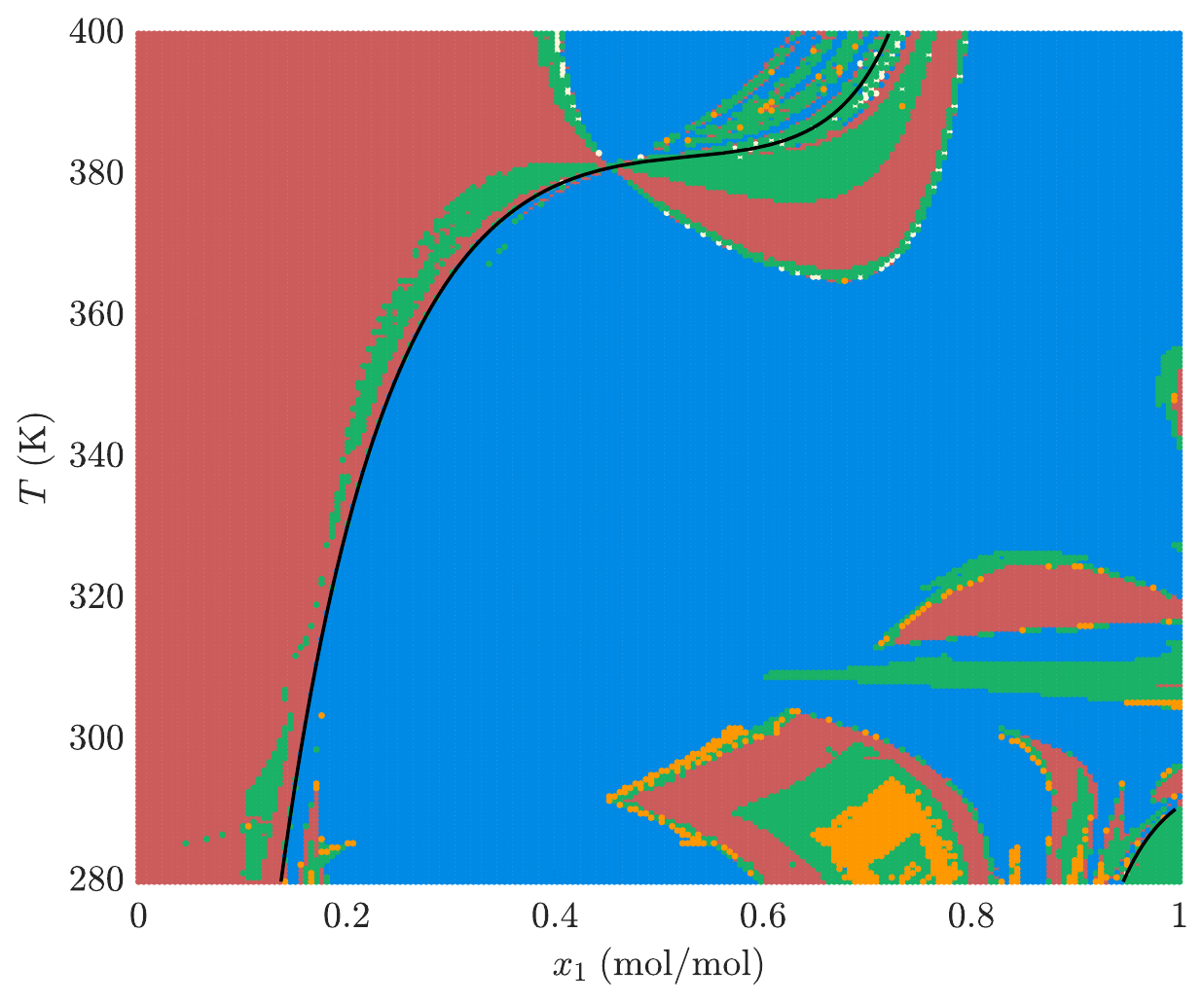}
}
  \caption{Global Newton method with error-oriented convergence criterion and adaptive trust region strategy (NLEQ-ERR) with 200 $\times$ 200 points. Red points: Azeotrope 1. Blue points: Azeotrope 2. Yellow points: singular Jacobian matrix. Green points: convergence to infinity. Orange points: oscillatory behavior. Black continuous line: critical curve.}
\label{fig:nleq_err}
\end{figure*}

The basin of attraction for the algorithm NLEQ-ERR is detailed in Figure \ref{fig:nleq_err}. Again, we observe a significant gain in terms of the non-convergence portions of the diagram, mainly close to the critical curve. On the other hand, a more complicated pattern is verified for molar fractions close to one.

\begin{figure*}[t]
 \centering
 \framebox[2.9in][c]{
 \includegraphics[width=2.7in]{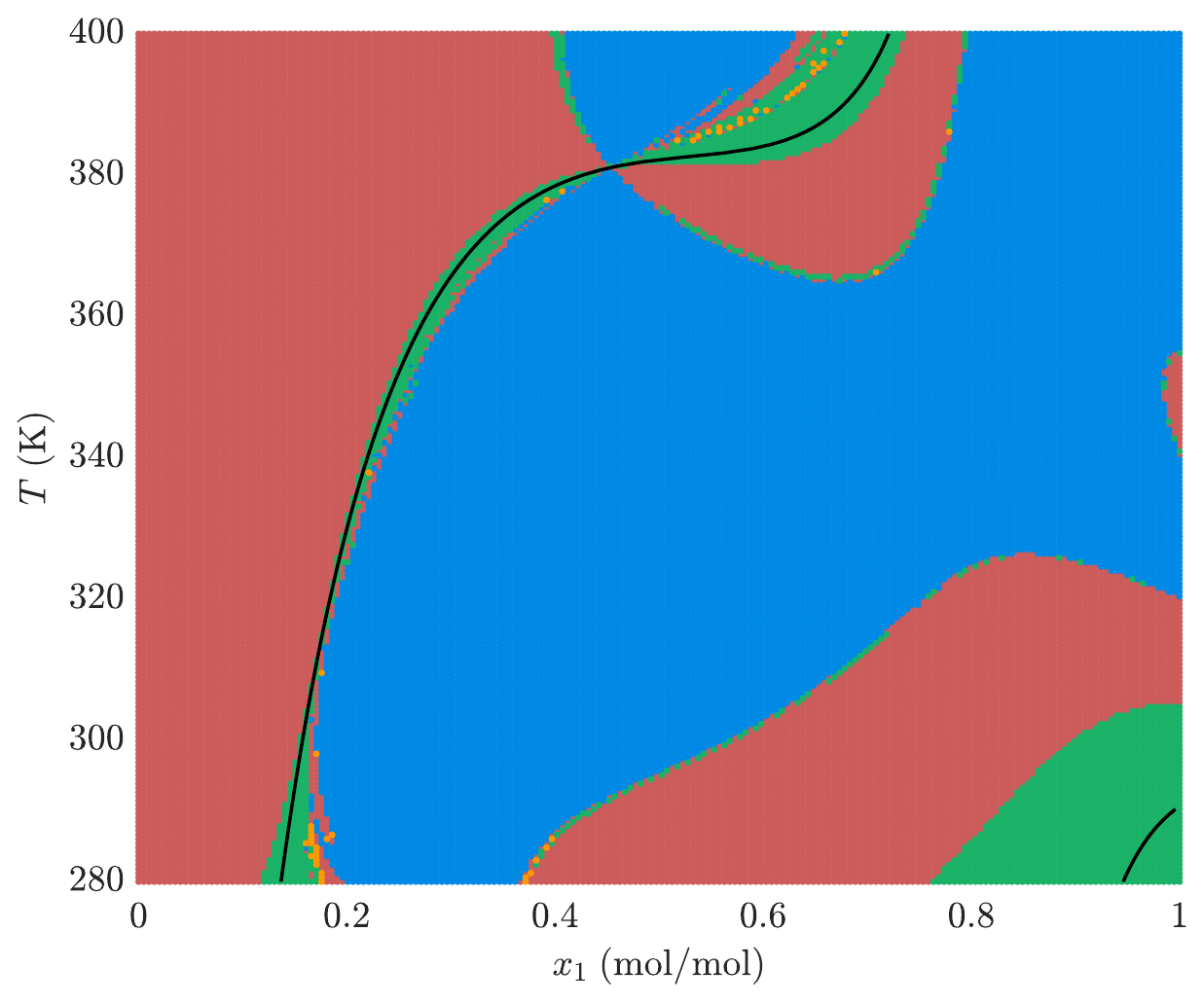}
}
  \caption{Jacobian-free Newton-Krylov method with 200 $\times$ 200 points. Red points: Azeotrope 1. Blue points: Azeotrope 2. Yellow points: singular Jacobian matrix. Green points: convergence to infinity. Orange points: oscillatory behavior. Black continuous line: critical curve.}
\label{fig:jfnk}
\end{figure*}

Figure \ref{fig:jfnk} illustrates the basin of attraction for the Newton-Krylov algorithm. As expected, there is no portion of the diagram corresponding to the singularity of the Jacobian. But there is a significant quantity of initial estimates in the neighborhood of the critical curve that converges to infinity (unwanted behavior, since we are mainly interested in the physical roots of the problem).

\begin{figure*}[t]
 \centering
 \framebox[2.9in][c]{
 \includegraphics[width=2.7in]{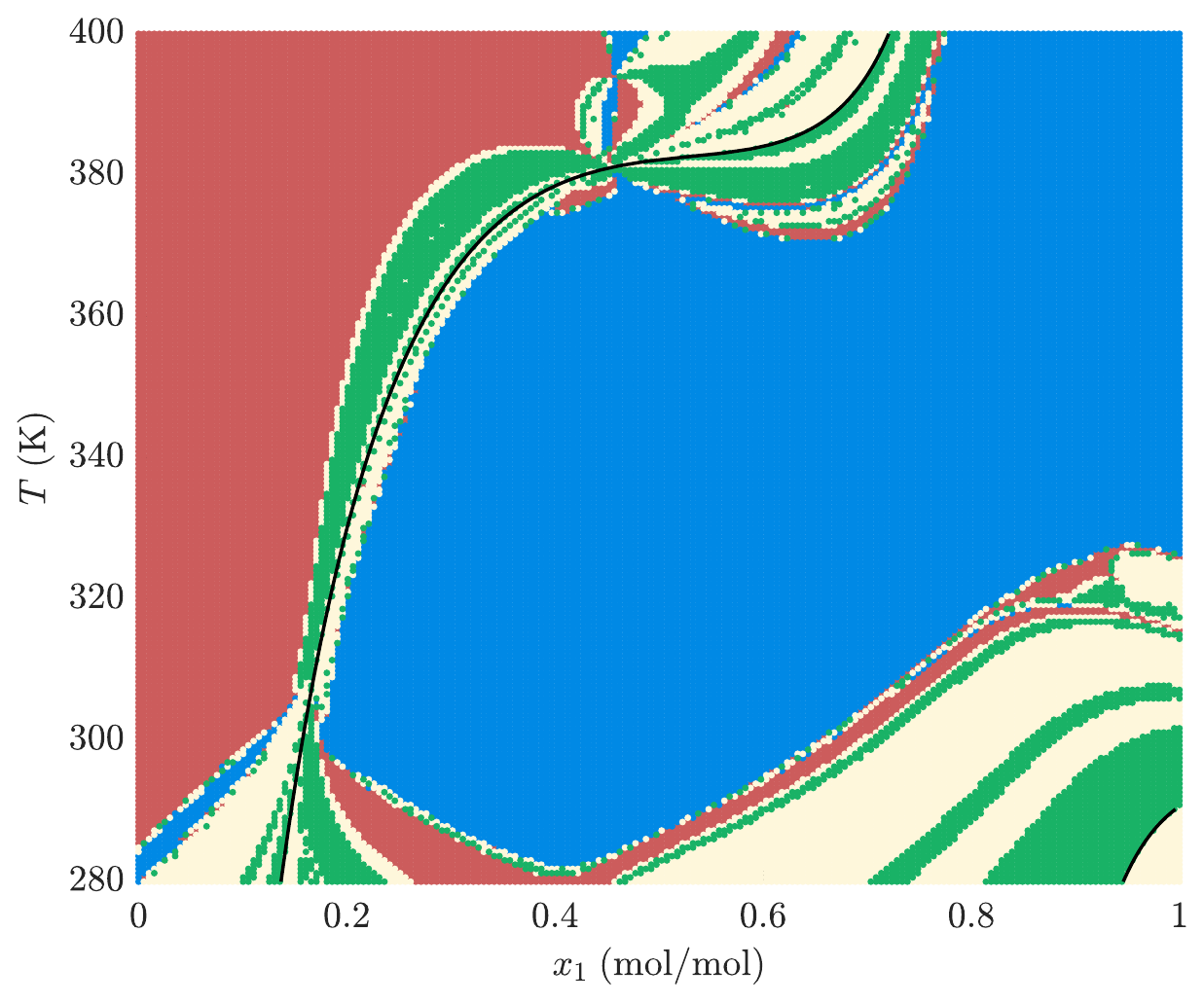}
}
  \caption{Method of \citet{Babajee2012} with 200 $\times$ 200 points. Red points: Azeotrope 1. Blue points: Azeotrope 2. Yellow points: singular Jacobian matrix. Green points: convergence to infinity. Orange points: oscillatory behavior. Black continuous line: critical curve.}
\label{fig:baba}
\end{figure*}

The results of the method of \citet{Babajee2012} are depicted in Figure \ref{fig:baba}. Clearly, the application of this algorithm was not recommended in this specific problem, since extremely complicated patterns were produced.

\begin{figure*}[t]
 \centering
 \framebox[2.9in][c]{
 \includegraphics[width=2.7in]{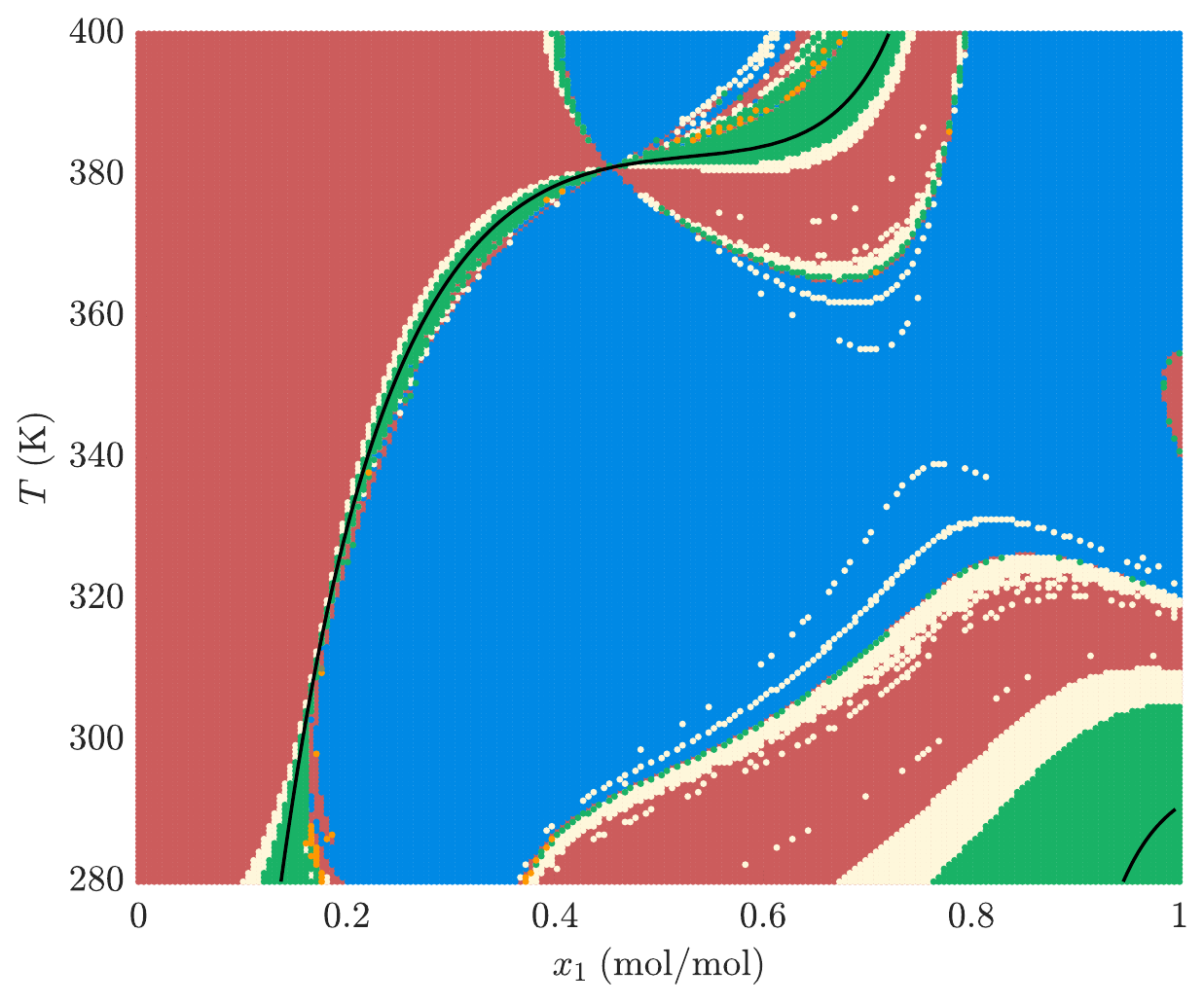}
}
  \caption{Method of \citet{Grau2011} with 200 $\times$ 200 points. Red points: Azeotrope 1. Blue points: Azeotrope 2. Yellow points: singular Jacobian matrix. Green points: convergence to infinity. Orange points: oscillatory behavior. Black continuous line: critical curve.}
\label{fig:grau}
\end{figure*}

\begin{figure*}[t]
 \centering
 \framebox[2.9in][c]{
 \includegraphics[width=2.7in]{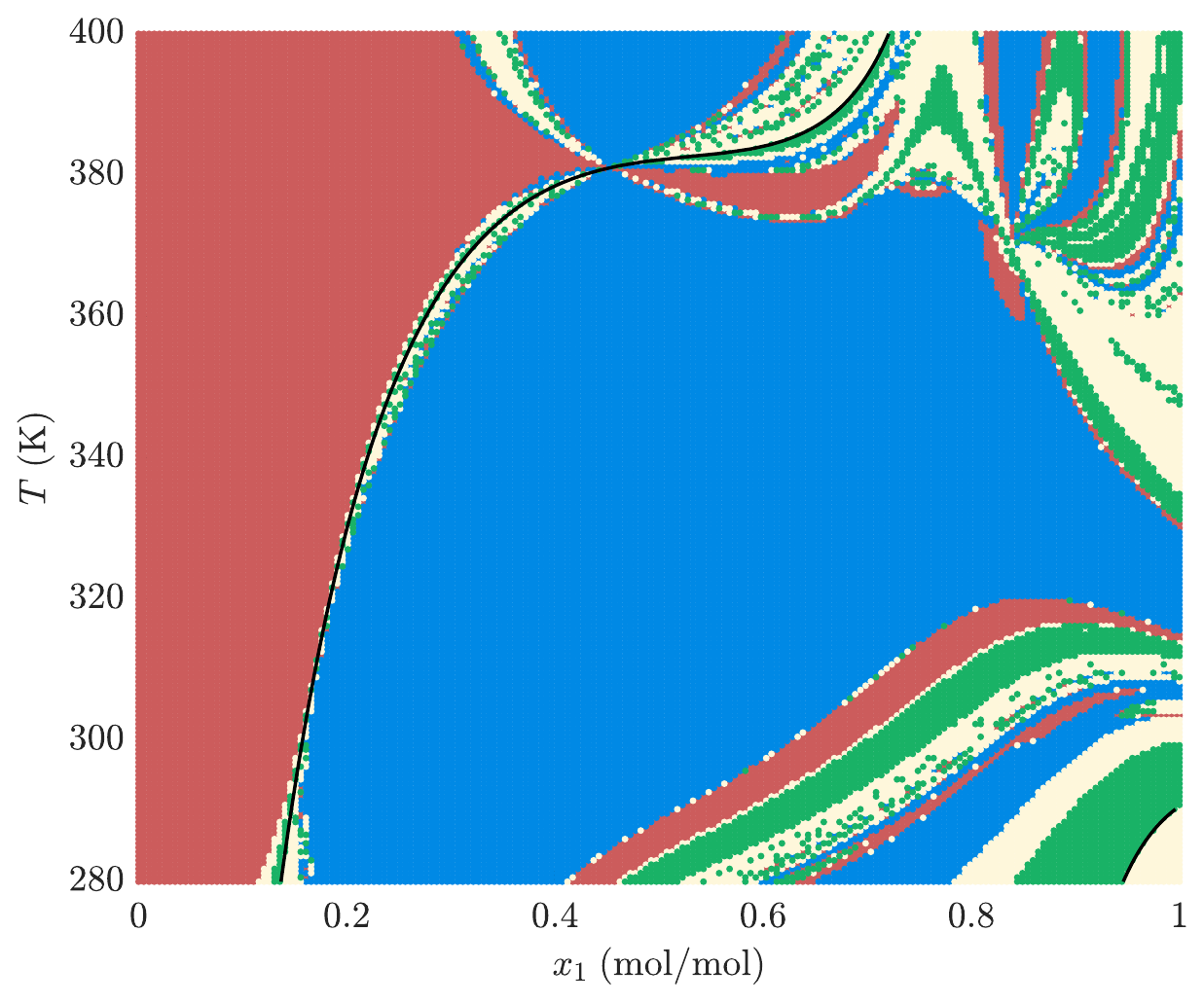}
}
  \caption{Method of \citet{Cordero2012} with 200 $\times$ 200 points. Red points: Azeotrope 1. Blue points: Azeotrope 2. Yellow points: singular Jacobian matrix. Green points: convergence to infinity. Orange points: oscillatory behavior. Black continuous line: critical curve.}
\label{fig:cordero}
\end{figure*}

An analysis of Figures \ref{fig:grau} and \ref{fig:cordero}, representing the basins of attraction for the methods of \citet{Grau2011} and \citet{Cordero2012}, respectively, indicates - qualitatively - a good behavior with respect of the non-convergence region in the proximity of the critical curves. Even with large portions of non-convergence/convergence to infinity close to the pure component 1 and a chaotic pattern for the algorithm of \citet{Cordero2012} at high temperatures, the vicinities of the critical curve are comparable to that obtained, for instance, with NLEQ-RES method.

\begin{figure*}[t]
 \centering
 \framebox[2.9in][c]{
 \includegraphics[width=2.7in]{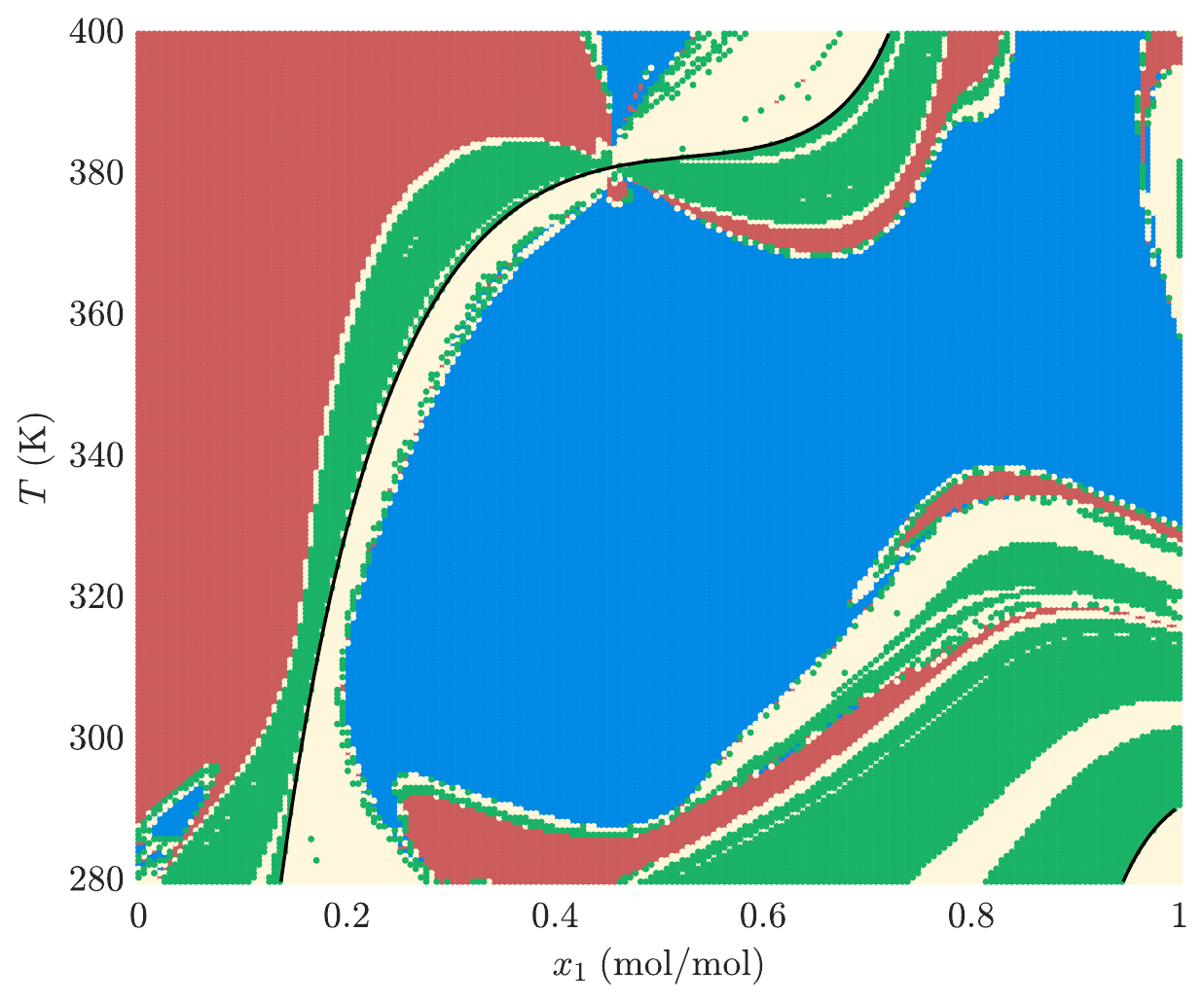}
}
  \caption{Fourth-order Newton-type method \citep{Madhu2017} with 200 $\times$ 200 points. Red points: Azeotrope 1. Blue points: Azeotrope 2. Yellow points: singular Jacobian matrix. Green points: convergence to infinity. Orange points: oscillatory behavior. Black continuous line: critical curve.}
\label{fig:fourth}
\end{figure*}

\begin{figure*}[t]
 \centering
 \framebox[2.9in][c]{
 \includegraphics[width=2.7in]{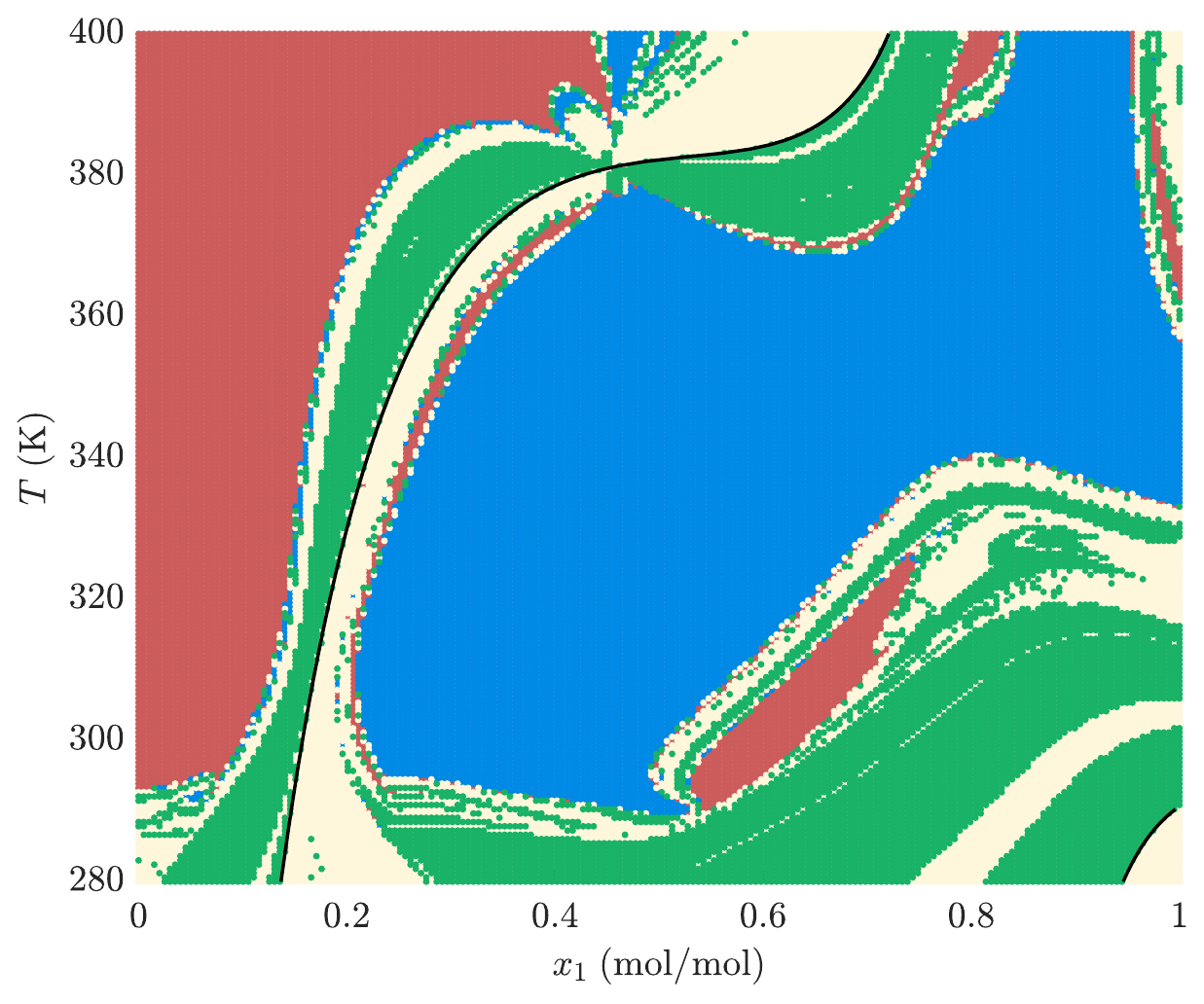}
}
  \caption{Fifth-order Newton-type method \citep{Madhu2017} with 200 $\times$ 200 points. Red points: Azeotrope 1. Blue points: Azeotrope 2. Yellow points: singular Jacobian matrix. Green points: convergence to infinity. Orange points: oscillatory behavior. Black continuous line: critical curve.}
\label{fig:fifth}
\end{figure*}

\begin{figure*}[t]
 \centering
 \framebox[2.9in][c]{
 \includegraphics[width=2.7in]{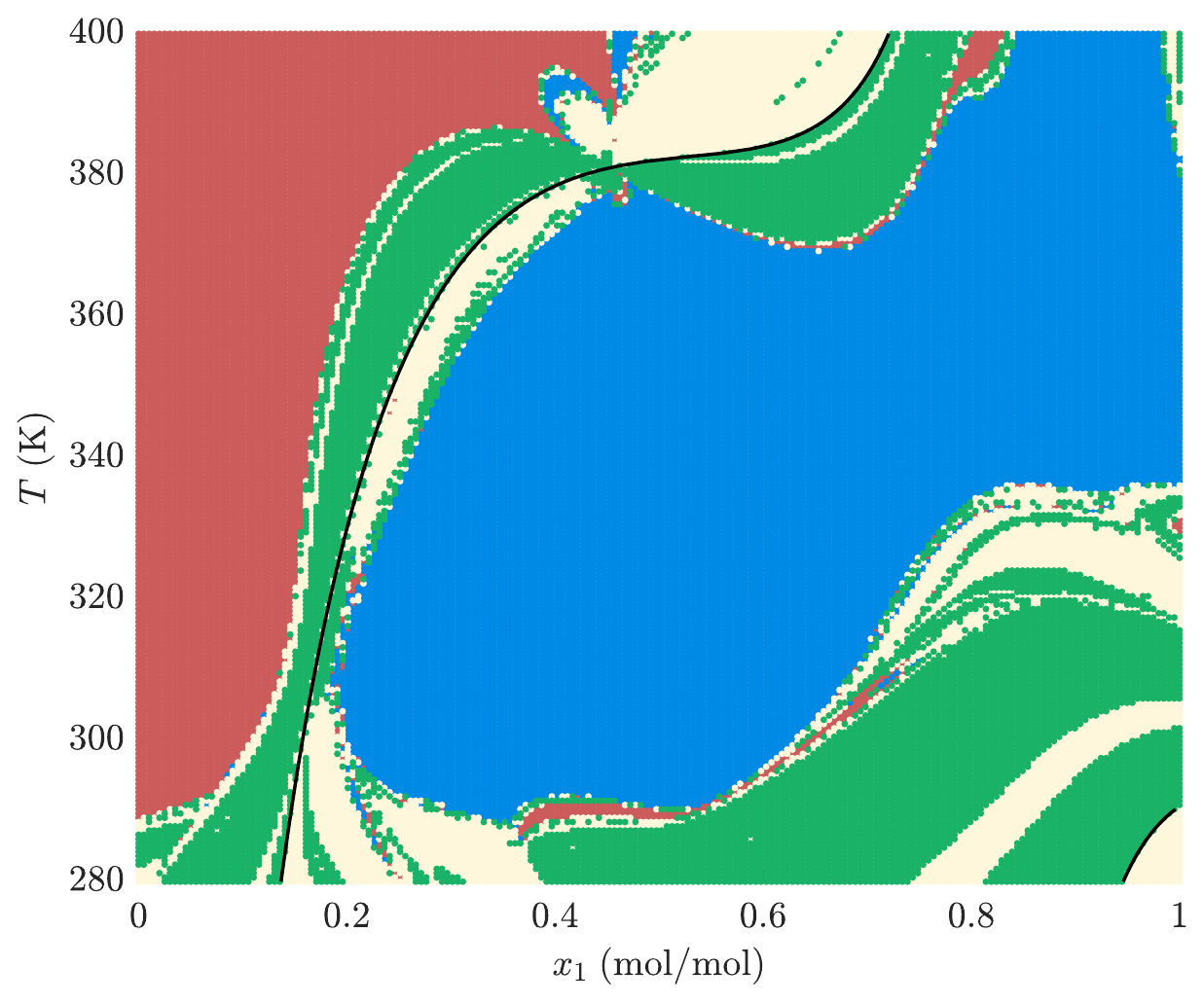}
}
  \caption{Sixth-order Newton-type method \citep{Madhu2017} with 200 $\times$ 200 points. Red points: Azeotrope 1. Blue points: Azeotrope 2. Yellow points: singular Jacobian matrix. Green points: convergence to infinity. Orange points: oscillatory behavior. Black continuous line: critical curve.}
\label{fig:sixth}
\end{figure*}

Figures \ref{fig:fourth} to \ref{fig:sixth} represent the basins of attraction for fourth to sixth order Newton's methods, respectively. We note a severe increase in the non-convergence region and the convergence to infinity - with comparison with the other methods - indicating an undesirable behavior for the identification of multiple solutions. \added{These results indicate that these algorithms only can be used when good initial estimates are available. Furthermore, the vicinity of the critical curve must be avoided, since we can note non-convergence regions in the neighborhood of the critical curve. Considering that, in practical problems, the explicit calculation of the critical curve is not performed, \textit{a priori}, the use of these expressions should be considered with care. Obviously, this recommendation cannot be extrapolated to other problems without a similar analysis; our main point here is to propose the application of the methodology for a particular class of problems.}

Basins of attraction illustrate -- in a qualitative way -- the convergence properties of the studied
methods, but we are also interested in the computational convergence order for the algorithms. Table \ref{tab:results} contains the results for the different algorithms considering four initial estimates. Two of them converged to Azeotrope 1 and the other ones converged to Azeotrope 2. We \replaced{analyzed}{analized} the number of function evaluations (NFE), the computational order of convergence $\rho_{avg}$ and the approximated computational order of convergence $\widehat {\rho }_{avg}$. The table also presents the convergence rate for each method, indicating the \replaced{percentage}{percentual quantity} of converged points in the 200 $\times$ 200 grid, considering only the physical solutions, i.e., without the convergence to
infinity.

The results displayed in Table \ref{tab:results} indicated that:

\begin{enumerate}
\item in a global point of view, CN and JFNK are the more ``robust'' methods, i.e., with lower rates of divergence;
\item the globalization techniques - represented in algorithms NLEQ-RES and NLEQ-ERR - provoked an
increase in NFE, but without effects in the convergence rate (in fact, we observed a small reduction in the convergence rate). This result is, in some sense, conflicting with that contained in the basins of attraction. In fact, we noted a reduction in the non-convergence condition for algorithms NLEQ-RES and NLEQ-EER only in the vicinities of the critical curve (mainly in NLEQ-ERR). On the other hand, somewhat complex patterns were verified in the other portions of the basin of attraction (for instance, close to the pure component 1 and low temperatures);
\item the approximated computational convergence order is severely affected by the initial estimate
considered in the results;
\item \added{the algorithms O5N e O6N, even presenting high computational convergence orders for a good initial estimate --- with $X^0 = (0.1,340)$ --- exhibited very low convergence rates in comparison with other techniques, and must be avoided in the absence of high quality initial estimates;}
\item qualitatively, higher values of $\rho_{avg}$ and $\widehat {\rho }_{avg}$ implied in lower convergence rates.
\end{enumerate}

\begin{table*}[h]
\centering
\resizebox{\textwidth}{!}{%
\begin{tabular}{l l l l l l l l l l l l l l}
\hline
\multirow{3}{*}{} & \multicolumn{6}{c}{} $\alpha = (0.0923864, 309.54)$ & \multicolumn{6}{c}{}   $\alpha = (0.2552517, 309.57)$  & \multirow{3}{*}{} \\ \cline{2-13} Method 
                  & \multicolumn{3}{c}{\replaced{${X^0} = (0.1,340)$}{${x_0} = (0.1,340)$}} & \multicolumn{3}{c}{\replaced{${X^0} = (0.2,380)$}{${x_0} = (0.2,380)$}} & \multicolumn{3}{c}{\replaced{$X^0 = (0.6,330)$}{$x_0 = (0.6,330)$}} & \multicolumn{3}{c}{\replaced{${X^0} = (0.8,360)$}{${x_0} = (0.8,360)$}} &     $C(\%)$               \\ \cline{2-13}
                  & NFE & $\rho_{avg}$     & $\widehat {\rho }_{avg}$     & NFE  & $\rho_{avg}$     & $\widehat {\rho }_{avg}$    & NFE  & $\rho_{avg}$     & $\widehat {\rho }_{avg}$    & NFE & $\rho_{avg}$     & $\widehat {\rho }_{avg}$     &                  \\ \hline
CN                & 5   & 1.7309 & 1.7388  & 7    & 1.8711 & 1.7803 & 4    & 1.5082 & 1.7272 & 5   & 1.2236 & 2.4567  & 0.8970             \\ \hline
NLEQ-RES          & 44  & 1.7309 & 1.7388  & 62   & 1.8711 & 1.7803 & 35   & 1.5082 & 1.7272 & 44  & 1.2236 & 2.4567  & 0.8881            \\ \hline
NLEQ-ERR          & 13  & 1.7312 & 1.7383  & 17   & 1.8712 & 1.7805 & 11   & 1.5092 & 1.7256 & 13  & 1.2219 & 2.4626  & 0.8881            \\ \hline
BA                & 16  & 0.9991 & 0.99    & 18   & 0.9991 & 0.9967 & 16   & 0.9995 & 0.9981 & 18  & 0.9995 & 0.9983  & 0.7548            \\ \hline
GS                & 9   & 1.3901 & 15.496  & 12   & 0.6486 & 3.6833 & 9    & 0.0033 & 3.981  & 9   & 1.5101 & 3.2736  & 0.8383            \\ \hline
CO                & 12  & 0.5813 & 17.259  & 12   & 3.0870 & 3.5727 & 12   & 0.6863 & 2.8384 & 12  & 2.1461 & 2.7981  & 0.8155            \\ \hline
O4N               & 8   & 1.5513 & 2.6390  & 10   & 1.0169 & 3.9754 & 8    & 0.0865 & 3.6147 & 8   & 1.2787 & 3.0926  & 0.6638            \\ \hline
O5N               & 9   & 3.2105 & 8.8649  & 12   & 2.4431 & 4.4258 & 9    & 1.2777 & 3.0538 & 12  & 0.0052 & 3.9642  & 0.5938            \\ \hline
O6N               & 9   & 2.3401 & 12.1644 & 12   & 1.4191 & 3.6947 & 9    & 0.8841 & 3.4398 & 9   & 2.5013 & 2.2356  & 0.6355            \\ \hline
JFNK              & 16  & 1.7301 & 1.7387  & 22   & 1.8703 & 1.7796 & 13   & 1.5084 & 1.7271 & 16  & 1.2223 & 2.46012 & 0.9315            \\ \hline
\end{tabular}}
\caption{Relationship among the number of function evaluations (NFE) and the convergence rate ($C(\%$)) of the analyzed methods, compared to the computational order of convergence ($\rho_{avg}$) and approximated computational order of convergence ($\widehat {\rho }_{avg}$), for four distinct initial estimates (\replaced{${X^0}$}{$x_0$}) and their respective solutions ($\alpha$).}
\label{tab:results}
\end{table*}

Considering the previous results -- in terms of basins of attraction, computational convergence orders
and convergence rates -- we will investigate, in a more detailed way, the algorithm NLEQ-ERR, considering two desirable properties of the method: (i) high convergence rates (close to the original Newton method and JFNK method); (ii) simple convergence patterns in the neighborhood of the critical curve.

\begin{figure*}[t]
 \centering
 \framebox[2.9in][c]{
 \includegraphics[width=2.7in]{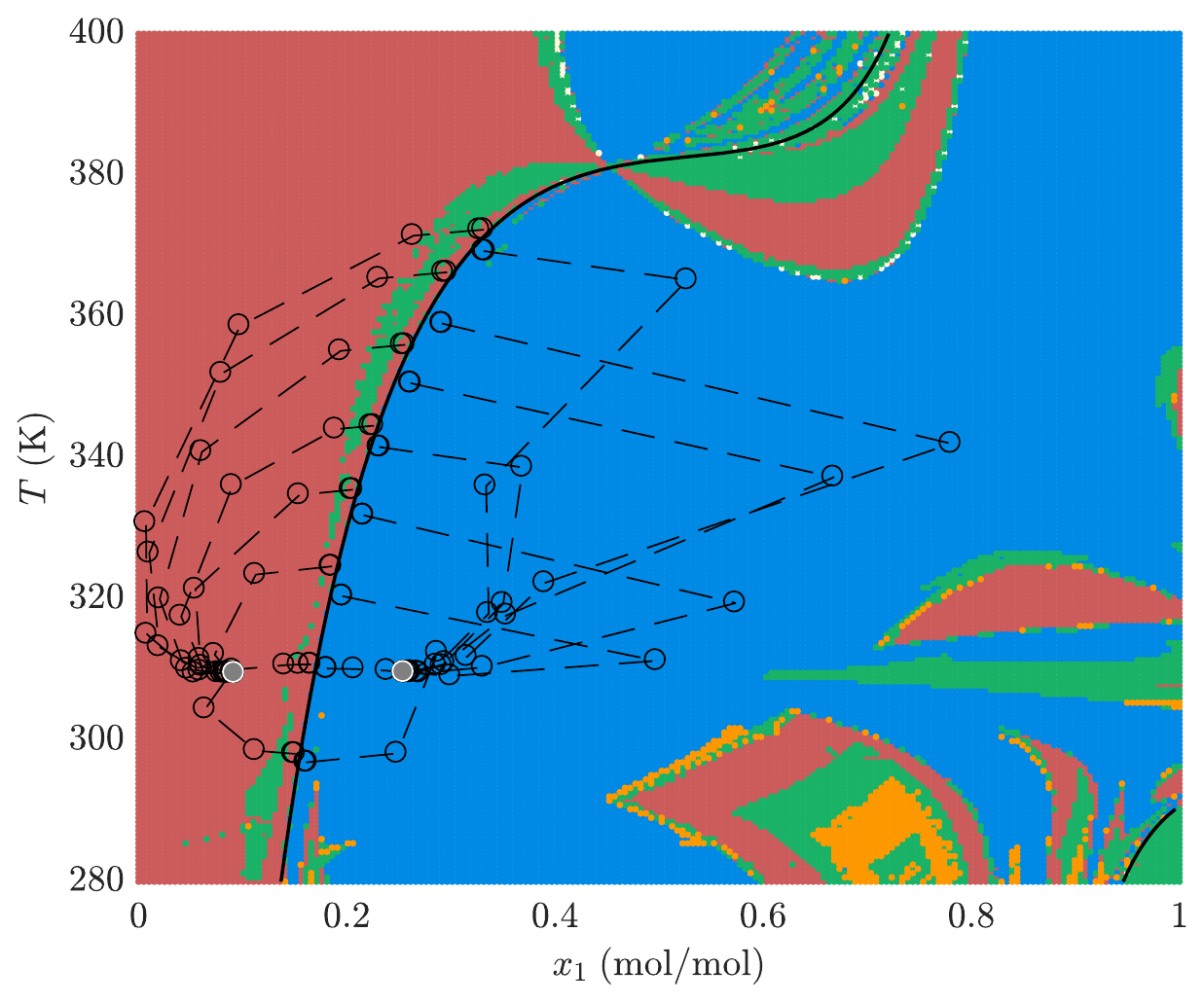}
}
  \caption{Iterative process of different runs of the NLEQ-ERR method, with initial estimates located in the vicinity of the critical curve and convergence for each solution of the double azeotropy problem, represented by the gray circles.}
\label{fig:newton_conv}
\end{figure*}

Figure \ref{fig:newton_conv} illustrates the iterative procedure - until convergence to the azeotropes, represented as gray circles - for some initial estimates in the vicinities of the critical curve for the algorithm NLEQ-ERR. In this situation we can note that the region with convergence to infinity (green points) is apart from the critical curve. Thus, with small perturbations in the critical curve we can obtain an adequate set of initial estimates. This behavior can be useful in the development of a new Newton-type algorithm, based in the identification of the critical curve of the problem.

\begin{figure*}[t]
 \centering
 \framebox[2.9in][c]{
 \includegraphics[width=2.7in]{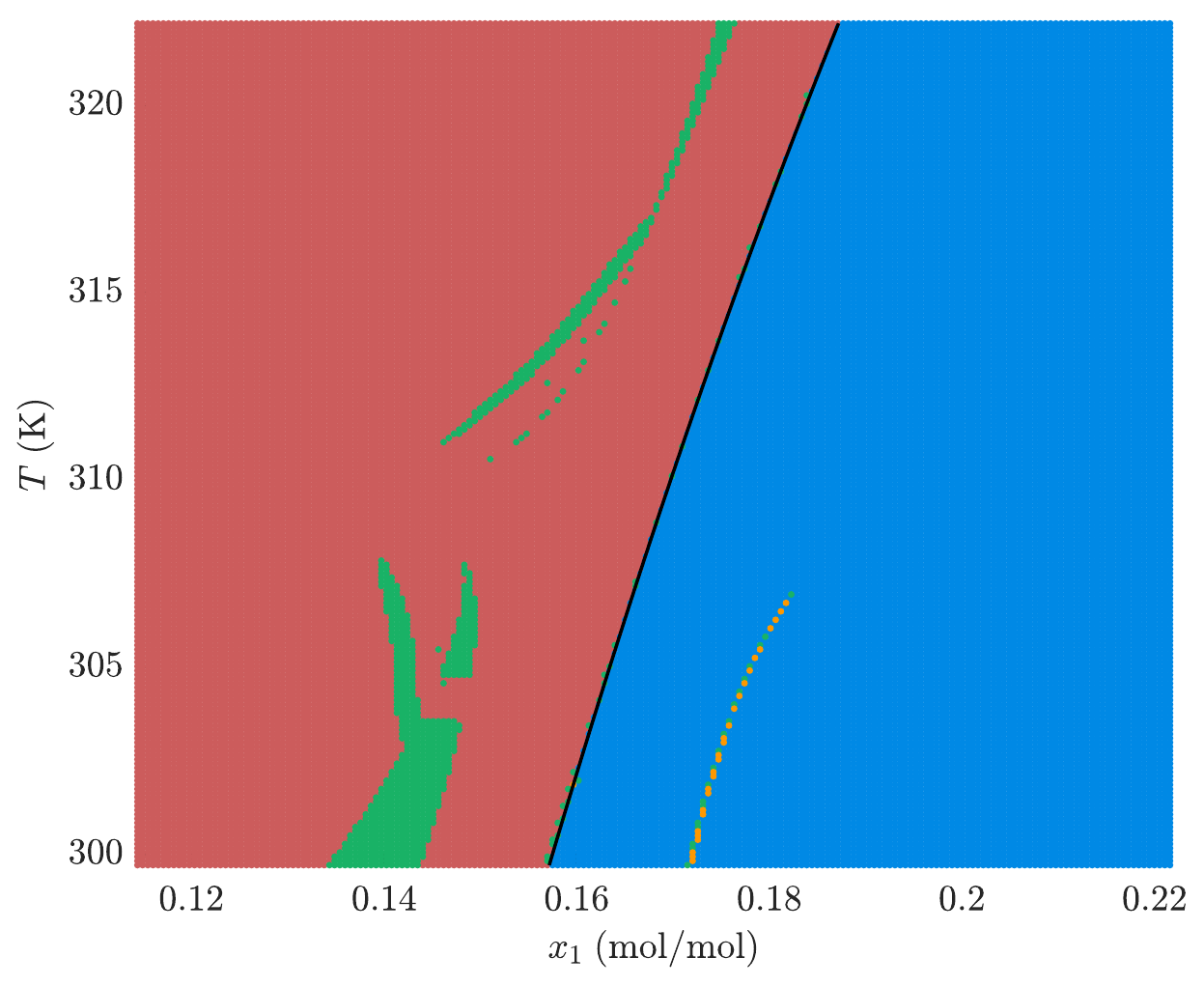}
}
  \caption{Amplification of the basin of attraction for global Newton method with error-oriented convergence criterion and adaptive trust region strategy (NLEQ-ERR) with 200 $\times$ 200 points. Red points: Azeotrope 1. Blue points: Azeotrope 2. Yellow points: singular Jacobian matrix. Green points: convergence to infinity. Orange points: oscillatory behavior. Black continuous line: critical curve.}
\label{fig:newton_amp}
\end{figure*}

Finally, Figure \ref{fig:newton_amp} presents an amplification of the region close to the critical curve for the method NLEQ-ERR. Notably, the identification of the critical curve -- in this problem -- permits the calculation, in a robust way, of the two physical solutions of the problem, fixing one variable and with small perturbations of the second one in opposite sides.

\section{Conclusions}
In this work we evaluated the basins of attraction for some Newton-type methods in the calculation of a double azeotrope (a problem in the plane), as well as presented the computational convergence orders for each algorithm.

The results indicated that:

\begin{enumerate}
\item basins of attraction and critical curves are closely related to the convergence properties of Newton-type methods;
\item some globalization techniques (mainly the algorithm NLEQ-ERR) promoted a more adequate convergence pattern in the vicinity of the critical curve, but without gains in terms of convergence rates;
\item high order methods implied -- in a qualitative analysis -- in more complicated behaviors in basins of attraction and lower convergence rates, \added{which suggests that these methods should only be used when good initial estimates --- far from the critical curve --- are available};
\item the identification of the critical curve, coupled with a well-designed algorithm (for a specific problem, such as NLEQ-ERR in this particular problem) can be used as a basis for the development of robust algorithms (in terms of capabilities for identification of all roots of the problem);
\item \added{this kind of analysis can be extended to other complex phase equilibrium problems with multiple solutions, serving as a guideline to the selection of robust root-finding algorithms.}
\end{enumerate}




\bibliographystyle{abbrvnat}
\bibliography{sample}

\end{document}